\documentclass[12pt]{article}

\def\C{{\bf C}}
\def\R{{\bf I\!R}}
\def\Z{{\bf Z\!\!Z}}
\def\Cl{{\mbox{\rm Cl}(V\oplus V^*)}}
\def\Vg{{\Lambda^{\cdot}(V)}}
\def\Vev{{\Lambda^{even}(V)}}

\def\r{{\rho}}
\def\dim{{\mbox{\rm \tiny dim}}}
\def\dimm{{\mbox{\rm  dim}}}
\def\QED{{\hspace*{3in} Q.E.D.}}
\def\At{{{\cal A}_{\theta}}}
\def\Atau{{{\cal A}_{\tau}}}
\def\n{{\nabla}}
\def\Lt{{L_{\theta}}}
\def\He{{{\cal H}eis}}
\def\Att{{{\cal A}_{\wt{\theta}}}}
 
\newtheorem{thm}{Theorem}[section]
\newtheorem{lem}{Lemma}[section]
\newtheorem{cor}{Corollary}[section]
\newtheorem{prop}{Proposition}[section]

\newtheorem{defi}{Definition}[section]

\def\wt{\widetilde}

\title{Projective modules
 over  non-commutative tori: classification of modules with constant
curvature connection. }

\author{Alexander Astashkevich\thanks{Renaissance Technologies, 
        E-mail address: ast@rentec.com}\,\,\,,
        Albert Schwarz\thanks{Department of Mathematics, 
        University of California, Davis, 
        E-mail address: schwarz@math.ucdavis.edu}  
\thanks{Research supported in part  by NSF  Grant No. DMS-9801009 and by
Mittag-Leffler Institute}}

\begin{document}

\maketitle

{\hspace*{2.4in} {To D. B. Fuchs on his 60-th birthday}\\

{\bf Abstract}.
We study finitely generated projective modules over noncommutative tori.
We prove that for every module $E$ with  constant curvature connection
the corresponding element $[E]$ of the K-group is a generalized quadratic
exponent and, conversely, for every positive generalized quadratic exponent
$\mu$ in the K-group one can find such a module $E$ with constant curvature
connection that $[E] = \mu $. In physical words we  give necessary and 
sufficient conditions for existence of 1/2 BPS states in terms of 
topological numbers.

\section{Introduction.}
\label{Sec:intro}

   In present paper we study projective modules over non-commutative tori.
(We always consider finitely generated
projective modules.)
Our main goal is to describe all modules that admit
constant curvature connections. It is well known that constant 
curvature connections correspond to  maximally supersymmetric BPS
fields \cite {cds} ; this means that we give conditions for existence of 
1/2 BPS states.

The main results of the paper are formulated in the following theorems.

\begin{thm}
\label{thm:M1}
Let $\At$ be a non-commutative torus.
Then for every projective $\At$ module $E$ with a constant
curvature connection corresponding element of the group  
$K_0(\At)$ is a generalized quadratic exponent.
Conversely,if $\mu$ is a positive generalized quadratic exponent in
$K_0(\At)$ then there exists such a projective module $E$ with constant
curvature connection that $[E]=\mu$.
(Here $[E]$ stands for the K-theory class of $E$.The definition of
generalized quadratic exponent will be given later. )
 
\end{thm}

\begin{thm}
\label{thm2}
Let $\At$ be an irrational non-commutative torus. In this case projective modules over $\At$ which admit 
constant curvature connection are in one-to-one
correspondence with positive generalized quadratic exponents in $K_0(\At)$.
\end{thm}

This theorem is an  immediate consequence of Theorem 1 and of the
following
 very strong result by M. Rieffel (see \cite{Rf1}):
 for irrational non-commutative torus $\At$  projective modules $E$ and $F$ 
are isomorphic if and only if the classes $[E],[F]\in K_0(\At)$ are
equal.

Our main results were formulated and partially proved in \cite {K-S},
Appendix D. 
It is assumed in \cite {K-S} that every linear combination of entries of
the matrix
$\theta$ is irrational. It is proved that in this case a projective module
can
be transformed into a free module by means of complete Morita equivalence
iff corresponding $K$-theory class is a generalized quadratic exponent. 
This statement can be used to prove
Theorem 1.1 in the conditions of Appendix D of \cite {K-S}).

The paper is organized as follows. In the introduction we remind
the main notions and results we need and explain how we plan to prove the 
main theorem.
In the section 2 we introduce the notion of generalized quadratic
exponent and we study its properties.
Section 3 is about integral generalized quadratic exponents
and finite dimensional representations of rational non-commutative tori.
In section 4 we present a proof of main results.

Let us remind the definition of a 
non-commutative torus (see \cite{R5} for more
details). Let $L$ be a lattice $\Z^n$ in a vector space $V^*=\R^n$.
Let $\theta$ be real valued skew-symmetric bilinear form on $R^n$.
We will think about $\theta$ as a two-form, that is an element
of $\Lambda^2 V$. Non-commutative torus $\At$ is the universal
$\C^*$-algebra generated by unitary operators $U_\alpha$, $\alpha\in L$
obeying relations
\begin{equation}
\label{eq:rel}
U_\alpha U_\beta=e^{\pi i \theta(\alpha, \beta)} U_{\alpha+\beta}.
\end{equation}
Any element from  $\At$ can be represented  uniquely by a sum  
$a=\sum_{{\bf \alpha}\in L} c_{\alpha}U_{\alpha}$, 
where $c_{\alpha}$ are complex numbers. Assigning to every $a\in\At$ the
coefficient $c_0$ in the representation above we obtain a canonical trace
$\tau$ on $\At$.

Let $\{e_i\}$ be a basis of $L$. One can say that $\At$ is the universal
$\C^*$-algebra generated by unitary operators $U_1,\cdots,U_n$
obeying the relations
\begin{equation}
\label{eq:relother}
U_i U_j=e^{2\pi i \theta(e_i, e_j)} U_j U_i.
\end{equation}
To check that these two  definitions are equivalent one should take  
$U_i=U_{e_i}$.

The transformations
 $\delta _k U_{e_k}=U_{e_k}$, $1\leq k\leq n$  ,
 $\delta _l U_{e_k}=0$, $k\not=l,1\leq k,l\leq n$ can be 
regarded as generators of Abelian Lie algebra $L_{\theta}$ of 
infinitesimal automorphisms on ${\cal A}_{\theta}$. 
We can naturally identify $\Lt$ with $V$.
Let us  remind the definition  of a connection in a ${\cal 
A}_{\theta}$-module following \cite{Con1} (we do not need the general
notion of connection \cite{Con2}). First we need the notion of a
smooth part of a  projective module.
 
Any element from $\At$ can be considered as a (generalized) function 
on the n-dimensional torus whose Fourier coefficients are $c_{\alpha}$  
(see above).  The space of smooth functions on $T^n$ forms a subalgebra 
of $\At$. We denote it by ${\At}^{smooth}$ and call it the smooth part of  
$\At$. If $E$ is a projective  $\At$ module one can  
define its smooth part $E^{smooth}$ in a similar manner (see \cite{Rf1}). 
A connection on  projective module $E$
can be defined as follows:
 
{\it $\At$ connection on a right $\At$ module $E$ 
is a linear map $\nabla :L_{\theta} 
\rightarrow \mbox{\rm End}_{\C}E$, 
satisfying the condition 
$$\n_{\delta}(ea)=(\n_{\delta}e)a+e(\delta(a)),$$ 
where $e\in E^{smooth}$, $a\in\At^{smooth}$, and $\delta\in L_{\theta}$. 
The curvature $F_{\mu,\nu}=[\n_{\mu},\n_{\nu}]$ of connection $\n$ is  
considered as a two-form on $\Lt$ with values in $\mbox{\rm End}_{\At}E$. 
(Here $\mbox{\rm End}_{\At}E$ stands for the space of 
endomorphisms of $\At^{smooth}$-module $E^{smooth}$ and $\mbox{\rm
End}_{\C}E$
 denotes the space of $\C$-linear
endomorphisms of $E^{smooth}$.)}
 
We always consider Hermitian modules and Hermitian connections.  
This means that if $E$ is a right $\At$ module it is equipped with  
$\At$ valued Hermitian inner product $\langle\cdot,\cdot\rangle$  
(for the detailed list of properties see \cite{Bl}) ;  all connections 
that we will consider should be compatible with this inner product. 

If $E$ is endowed with a $\At$-connection, then one can define a 
Chern character
\begin{equation}
\label{eqchern}
\mbox{\rm ch}(E)=\sum_{k=0}\frac{\hat{\tau}(F^k)}{(2\pi i)^k k!}=
\hat{\tau}(e^{\frac{F}{2\pi i}}),
\end{equation}
where $F$ is a curvature of a connection on $E$, and $\hat{\tau}$ is the
canonical trace on $\hat{A}=\mbox{\rm End}_{\At}(E)$ 
(we use that
 $\At$ is equipped with a canonical trace $\tau=c_0$). One can consider
 $\mbox{\rm ch}(E)$ as an element in the  Grassmann algebra 
$\Lambda^\cdot(\Lt^*)=\Lambda^\cdot (V^*)$. 
We have a lattice $L$ in $V^*$. Thus we can talk about integral
elements in $\Lambda^\cdot (V^*)$ which are just the elements of
$\Lambda^\cdot L$.
In the commutative case ${\rm ch}(E)$ is integral.
In non-commutative case this is wrong, but there exists an integral
element $\mu(E)\in \Lambda^\cdot (V^*)$ related to 
 $\mbox{\rm ch}(E)$  by the 
 formula (see \cite{Ell}, \cite{Rf1})
\begin{equation}
\label{eqrelchmu}
\mbox{\rm ch}(E)=e^{\iota(\theta)}\mu(E),
\end{equation}
Here ${\iota(\theta)}$ stands for the operation of contraction with 
${\theta}$ considered as an element of $\Lambda^2 V$.
In particular, formula (\ref{eqrelchmu}) means that
$e^{-\iota(\theta)}\mbox{\rm ch}(E)$ is an integral element 
of $\Lambda^\cdot (V^*)$. The group $\Lambda^{even} L$ can be naturally
identified with the group $K_0(\At)$. Moreover
$\mu(E)$ is the class of the module $E$ in the $K_0(\At)$ group (see
\cite{Ell}).

Let us remind that the element $\mu\in K_0(\At)$ is called positive
if $(e^{\iota(\theta)}\mu)_{(0)}>0$ (the zero component is positive).
A well known theorem of M. Rieffel (see \cite{Rf1}) says that if $\theta$
is irrational then every positive element of $\mu$ is represented by a
 projective module over $\At$.

Let $E$ be a  projective $\At$ module with a 
constant curvature connection $\n$. Denote by $F$ the curvature of $\n$.
Then since $F$ is a 2-form with values in $\C$ we obtain that
$\hat{\tau}(F^k)=\hat{\tau}(1) F^k$. The number $\hat{\tau}(1)$ is
called the dimension of the module $E$ and we denote it by $d_E$.
Then the formula \ref{eqchern} becomes
\begin{equation}
\label{eq:ch}
ch(E)=d_E e^{\frac{F}{2\pi i}}.
\end{equation}
We see that in this case $ch(E)$ is a quadratic exponent (i.e. an
expression of the form $Ce^a$ where $C$ is a constant and $a\in \Lambda
^2(V^+)$). It follows from (4) and from this fact that $\mu (E)$ is a
generalized quadratic exponent (i.e. a limit of quadratic exponents). This
gives a proof of the first statement of Theorem \ref{thm:M1}. The proof of
the second statement of this theorem is based on the study of 
 generalized quadratic
exponents in sections 2 and 3. In section 3 we will
study integral generalized quadratic exponents keeping in mind that
$K_0(\At)$ is exactly 
the integral lattice in $\Vev$, where $V=\Lt^*$. We will prove some
auxiliary technical results saying that something is rational
or integral which we will use in our construction of the module
in section 4.\\

In section 4 we will construct a desired module together
with constant curvature connection in four steps.

First we find explicitly the curvature $F$ as a 2-form on $\Lt$.

Second, we construct some spaces of functions on $\R^p\times \Z^{n-2p}$
with action of  generators of 
some non-commutative torus ${\cal A}_{\wt{\theta}}$ and  with  constant
curvature connection having the curvature form $F$.
We do not construct an ${\cal A}_{\wt{\theta}}$ module at this step.
Moreover, we even will not specify what space of functions we will
take. 

Third, we will check that our construction in the previous step
is just a particular case of Rieffel's construction in \cite{Rf1}
where he constructs  ${\cal A}_{\wt{\theta}}$  
projective modules. So we can construct the desired module $\wt{E}$ over
${\cal A}_{\wt{\theta}}$ using  Rieffel's construction.

Fourth, we will see that $\tau=\theta-\wt{\theta}$ is a rational element
of $\Lambda^2 \Lt=\Lambda^2 V^*$. 
Also, we will use Rieffel's  explicit calculation
 of $[\wt{E}]$ (of the class of $\wt{E}$ in 
$K_0({\cal A}_{\wt{\theta}})\subset \Vev$ )to find a simple relation
between $[\wt{E}]$ and $\mu$.
Finally, we show that we can construct
a  projective module $E$ with constant curvature 
connection over $\At$ such that $[E]=\mu$ by taking $E$ to be a tensor
product of $\wt{E}$ by some finite dimensional module $M$ over $\Atau$.

\section{Generalized quadratic exponents.}
\label{sec:gen}

In this section we introduce  generalized quadratic exponents
and study their properties.

Let $V$ be a finite dimensional vector space over $\R$. Let $V^*$ be
 a dual space. Then the space $V\oplus V^*$ has a natural symmetric
bilinear product given by
\[
\langle (x_1,y_1),(x_2,y_2)\rangle=y_2(x_1)+y_1(x_2),
\]
where $x_1,x_2\in V$ and $y_1,y_2\in V^*$. Consider a Clifford algebra
$\Cl$. It naturally acts on the vector space $\Vg$; we denote this action
by $\r$. 
Note, that there is a natural inclusion $i$ of $V\oplus V^*$ into $\Cl$.

\begin{defi}
\label{defi:gqe}
An element $q\in \Vg$ is called a generalized quadratic exponent
if there exists a maximal isotropic  subspace $U\subset V\oplus V^*$
such that for any $x\in U$ we have $\r(x)q=0$.
\end{defi}

If the projection of $U$ onto $V^*$ is bijective we can represent $U$ as a
graph of a linear operator $a: V^*\rightarrow V$. The operator $a$ is
antisymmetric; it can be considered as an element of $\Lambda ^2(V)$. The
element $q$ can be represented in the form const$\cdot e^a$, i.e. it is a
quadratic exponent. The set of maximal isotropic subspaces we just
considered is
dense in the set of all maximal isotropic subspaces; this means that
quadratic exponents are dense in the set of all generalized quadratic
exponents. 

In the next proposition we will describe all possible generalized quadratic
exponents.

\begin{prop}
\label{prop:gqe}
Let $q\in\Vg$ be a generalized quadratic exponent. Then there exists
a subspace $W\subset V$, a non-degenerate element $\wt{q}_1\in \Lambda^2(V/W)$
and nonzero element $w\in \Lambda^{\dim W}W$ such that 
\[
q=w\wedge e^{q_1}.
 \]
where $q_1\in \Lambda^2(V)$ is any preimage of $\wt{q}_1\in \Lambda^2(V/W)$
under the natural projection from $\Lambda^2(V)$ onto $\Lambda^2(V/W)$.
\end{prop}

{\bf Proof:} Let $U$ be the maximal isotropic subspace corresponding to $q$.
It is easy to see that we can choose a basis $\{\xi_i\}$ 
of $V$ and a dual basis
$\{\eta_i\}$ of $V^*$ such that $U$ is spanned by the vectors
$\eta_1-\sum_{i=1}^j a_{1,i} \xi_i,...,
\eta_j-\sum_{i=1}^{j} a_{j,i} \xi_i,
\xi_{j+1},...,\xi_{\dim W}$. 
Thus, $q$ satisfies the following system of equations:
\[
\begin{array}{l}
\frac{\partial q}{\partial \xi_1}-(\sum_{i=1}^j a_{1,i} \xi_i)\wedge q=0\\
\cdot\\
\cdot\\
\cdot\\
\frac{\partial q}{\partial \xi_j}-(\sum_{i=1}^j a_{j,i} \xi_i)\wedge
q=0\\[12pt]
\xi_{j+1}\wedge q=0\\
\cdot\\
\cdot\\
\cdot\\
\xi_{\dim W}\wedge q=0
\end{array}
\]
The partial derivatives in this system are understood as left derivatives
in the sense of superalgebra.
 
It is easy to see that any solution of this system is of the form
\[
C\cdot\xi_{j+1}\wedge\cdots\wedge\xi_{\dim W}\wedge
e^{\sum_{k=1}^j\sum_{l=1}^j a_{k,l} \xi_k\wedge \xi_l},
\]
where C is a constant. The proposition follows easily from the
above formula. $W$ is the subspace spanned by $\xi_{j+1},\cdots,\xi_{\dim
W}$ and 
$\wt{q}_1$ is the projection of 
${\sum_{k=1}^j\sum_{l=1}^j a_{k,l} \xi_k\wedge \xi_l}$.\\
\QED

$\Vg$ is a graded vector space. If $q\in\Vg$ let us denote by
$q_{(i)}\in \Lambda^i V$ the projection $q$ on $\Lambda^i V$.

\begin{cor}
\label{cor:const}
Let $q$ be a generalized quadratic exponent. If $q_{(0)}$ is 
not zero then there is a non-degenerate element $a\in\Lambda^2 V$
and a nonzero real number $C$ such that 
\[
q=C e^{a}.
\]
\end{cor}

{\bf Proof:} Immediately follows from Proposition \ref{prop:gqe}.\\
\QED

Let $b\in \Lambda^2 (V^*)$. Then $b$ acts naturally on $\Vg$. If we choose
basis  $\{\xi_i\}$ in $V$ then we can write the action of $b$ as
$\sum_{k,l} b_{k,l}\frac{\partial}{\partial \xi_k \partial \xi_l}$. 
Another way of thinking is to think about $b$ as an element of $\Cl$.
We have a canonical map from $\Lambda^\cdot (V^*)$ to $\Cl$ since
$V^*$ is isotropic subspace in $V\oplus V^*$. Then the action
of $b$ is simply $\r(b)$.

\begin{prop}
\label{prop:closure}
Let $q$ be a generalized quadratic exponent and $b$ any element
in $\Lambda^2 (V^*)$. Then $e^{\r(b)}q$ is a generalized quadratic
exponent.
\end{prop}

{\bf Proof:}
We will reduce the proposition to the case where $b$ is decomposable. Since, 
$\r(b)=\sum_{k,l} b_{k,l}\frac{\partial}{\partial \xi_k \partial \xi_l}$
in some basis $\{\xi_i\}$ of $V$ and the operators 
$\frac{\partial}{\partial \xi_k \partial \xi_l}$ commute it is enough
to prove the proposition in the case when 
$\r(b)=c\frac{\partial}{\partial \xi_k \partial \xi_l}$, where c is a real
number. In this case
$e^{\r(b)}=1+\r(b)$. 

Our goal is to show that there exists a subspace $\wt{W}\subset V\oplus V^*$
such that for any $x\in\wt{W}$ we have $\r(x)(q+\r(b)(q))=0$.

Let $W$ be a subspace of $V\oplus V^*$ such that if $x\in W$ then $\r(x)q=0$.
We can choose a basis $\{v_1+w_1,\cdots,v_k+w_k,v_{k+1},\cdots,v_{\dim W}\}$
of $W$, where $v_i\in V^*$ and $w_i\in V$, and the vectors $\{w_i\}$ are
linearly independent. 

A simple calculation shows that
\[
\r(v_l)(q+\r(b)(q))=0+\r(v_l)\r(b)(q)=\r(b)\r(v_l)(q)=0,
\]
for $l>k$. Also, we can easily see that for $l<k+1$
\[
\begin{array}{l}
\r(v_l+w_l)(q+\r(b)(q))=\r(v_l+w_l)\r(b)(q)=\\[12pt]
[\r(v_l+w_l),\r(b)](q)+\r(b)
\r(v_l+w_l)(q)
=[\r(w_l),\r(b)](q)=\r(\iota(w_l)b)q,
\end{array}
\]
where $\iota(w_l)$ is plugging the vector $w_l$ in the 2-form $b$.
$\iota(w_l)b$ is an element of $V^*$.  Since $b^2=0$ we see that
$0=\iota(w_l)(b^2)=2(\iota(w_l)b)b$. Thus,
\[
\r(\iota(w_l)b)q=\r(\iota(w_l)b+(\iota(w_l)b)b)q=\r(\iota(w_l)b)(q+\r(q)).
\]
Therefore, we see that
\[
\r([v_l-\iota(w_l)b]+w_l)(q+\r(b)q)=0.
\]
Denote by $\wt{W}$ the subspace of $V\oplus V^*$ spanned by
the vectors
$[v_1-\iota(w_1)b]+w_1,\cdots,[v_k-\iota(w_k)b]+w_k,v_{k+1},...,
v_{\dim W}$. It is easy to check that $\wt{W}$ is a maximal isotropic
subspace of $V\oplus V^*$ and we showed that $\r(x)(q+\r(b)q)=0$ 
for any $x\in \wt{W}$. Thus, $q+\r(b)q$ is a generalized quadratic exponent.\\
\QED



\section{Integral generalized quadratic exponents.}
\label{sec:int}

In this section we study integral generalized
quadratic exponents and we prove a couple of auxiliary propositions 
that we will use in our construction.

Let $V$ be a finite dimensional vector space and let $L$ be a lattice
in it. Denote by $n$ the dimension $V$. 
Then, $V\cong \R^n$ and $L\cong \Z^n$.  We denote the dual lattice to $L$
by $L^*$. Obviously $L^*\subset V^*$. We call an element of $\Vg$
integral if it lies in $\Lambda^{\cdot}L$.

Define  $U_\mu$  a subspace
of $V^*$ as follows:
\[
U_\mu=\left\{ x\in V^*\,\mid \, \iota(x)\mu=0\right\}.
\]
Denote by $W_\mu\subseteq V$ the orthogonal complement to $U_\mu$.

\begin{prop}
\label{prop:intgqe}
Let $\mu\in \Vg$ be an integral generalized quadratic exponent.
Then $L_\mu=L\cap W_\mu$ is a lattice in $W_\mu$. We can identify
$\Lambda ^{\dim W_\mu} L_\mu$ with $\Z$ (the isomorphism is not canonical
but it is specified up to a sign). Let $\alpha\in \Lambda ^{\dim W_\mu}
L_\mu$ be 
a volume form (an
element that corresponds to $1$ under the isomorphism with $\Z$).
Then $\mu_{(\dim W_\mu)}=N \alpha$, where $N$ is a nonzero integer.
\end{prop}
{\bf Proof:} Let $\mu\in \Vg$ be an integral generalized quadratic exponent.
Let $k$ be the largest integer such that $\mu_{(k)}\neq 0$ and
for all $l>k$ we have $\mu_{(l)}=0$. 
{}From Proposition \ref{prop:gqe} easily follows that 
$U_\mu=\left\{ x\in V^*\,\mid \, \iota(x)\mu=0\right\}=
\left\{ x\in V^*\,\mid \, \iota(x)\mu_{(k)}=0\right\}$.
Moreover, from Proposition \ref{prop:gqe} follows that $\mu_{(k)}$
is a decomposable element of $\Vg$ and that
$k=\dimm V-\dimm U_\mu=\dimm W_\mu$. Since $\mu$ is integral $\mu_{(k)}$
is also integral. Thus, the subspace $U_\mu$ is spanned by $U_\mu\cap L^*$.
Therefore, $U_\mu\cap L^*$ is a lattice in $U_\mu$. This immediately
implies that $L_\mu=W_\mu\cap L$ is a lattice in $W_\mu$ since $W_\mu$ is
the orthogonal complement to $U_\mu$. 

{}From the above discussion it is easy to see that 
$\Lambda ^{\dim W_\mu} (W_\mu\cap L)=\Lambda ^{\dim W_\mu} L_\mu\cong \Z$.
Since, by the definition $\alpha$ corresponds to $\pm 1$ under such an 
isomorphism and $\mu_{(\dim W_\mu)}=\mu_{(k)}$ is an integral element
we obtain that $\mu_{(\dim W_\mu)}=N \alpha$ for some integer $N$.\\
\QED

Under the conditions in the above proposition we can easily find
a complement $\wt{L}_\mu$ ($\wt{L}_\mu\cong \Z^{n-\dim W_\mu}$)
to $L_\mu$ in $L$. It is not unique but we do 
not care about that. Let $Y_\mu$ be the subspace of $V$ spanned by 
$\wt{L}_\mu$. Then it is obvious that $V=W_\mu\oplus Y_\mu$ and
$L=L_\mu\oplus \wt{L}_\mu$.

Next results will be used in the construction in section
\ref{sec:constr}. Since, they do not use any theory of
non-commutative tori we state them here. But they need some explanation
concerning their origin. 

Let $\mu\in\Vev$ be an integral generalized
quadratic exponent which will be an element of $K_0$ representing
a projective module of $\At$. 
We can think about $\theta$ as an element of $\Lambda^2 V^*$. 
If there exists a projective module $E$  over $\At$ with constant 
curvature 
connection such that $[E]=\mu$ then by the result of G. Elliott (see
\cite{Ell})
$d e^{\frac{F}{2\pi i}}=e^{\iota(\theta)}\mu$, where $d$ is the dimension
of the module.  $\wt{\theta}$ which satisfies the conditions
of the lemma below will be constructed in section \ref{sec:constr}.

\begin{lem}
\label{lem:ration}
Let $\mu\in\Vev$ be an integral generalized
quadratic exponent.
Let us assume that we fixed the isomorphism between 
$\Lambda ^{\dim W_\mu} L_\mu$ and  $\Z$  (see Proposition \ref{prop:intgqe})
so that 
$\mu_{(\dim W_\mu)}=N\alpha$ with $N$ being natural number (here
$\alpha\in \Lambda ^{\dim W_\mu} L_\mu$ corresponds to $1$ in $\Z$). 
Let $\theta$ and $\wt{\theta}$ be elements of $\Lambda^2 V^*$ such that
$\theta-\wt{\theta}$ is zero on $V\otimes Y_\mu$ (that is if $X_\mu\subset V^*$
is the orthogonal complement to $Y_\mu$ then 
$\theta-\wt{\theta}\in \Lambda^2 X_\mu$). Assume that
\begin{equation}
\label{eq1}
e^{\iota(\theta)}\mu=c e^{\iota(\wt{\theta})}\alpha,
\end{equation}
where $c$ is a real number. Then $c=N$, that is
\begin{equation}
\label{eq2}
e^{\iota(\theta)}\mu=N e^{\iota(\wt{\theta})}\alpha,
\end{equation}
and $N (\theta-\wt{\theta})$ is an integral element of $\Lambda^2 X_\mu$.
\end{lem}

{\bf Proof:} Let us denote $\dimm W_\mu$ by $k$. Then formula
\ref{eq1} implies that 
\[
\mu_{(k)}=(e^{\iota(\theta)}\mu)_{(k)}=c (e^{\iota(\wt{\theta})}\alpha)_{(k)}=
c \alpha_{(k)}.
\]
Thus, $c=N$ and we proved formula \ref{eq2}. From formula \ref{eq2}
it easily follows that
\begin{equation}
\label{eq3}
\mu=N e^{\iota(\wt{\theta}-\theta)}\alpha.
\end{equation}
This means that $\mu_{k-2}=N \iota(\wt{\theta}-\theta)\alpha$.
$\mu_{k-2}$ is an integral element.
$\wt{\theta}-\theta$ is in $\Lambda^2 X_\mu$ which is dual to $W_\mu$
and $\alpha$ is a nonzero element of $\Lambda^{\dim W_\mu}W_\mu$.
Thus, $N(\wt{\theta}-\theta)$ is an integral element of $\Lambda^2 X_\mu$.\\
\QED

\subsection{Modules over the rational non-commutative tori.}
\label{subsec:ration}

Let us assume that conditions of  Lemma \ref{lem:ration} are satisfied.
We denote ${\theta}-\wt{\theta}$ by $\tau$. Let us remind the definition
of $\Atau$. $\Atau$ is a universal $\C^*$ algebra having unitary 
generators $U_\beta$, where $\beta\in L$, obeying the relations
\[
U_{\beta_1}U_{\beta_2}=e^{\pi i \tau(\beta_1,\beta_2)} U_{\beta_1 +\beta_2}.
\]
We can reformulate this definition in a slightly different way.
Let $\beta_1,\cdots,\beta_n$ (where $n=\dimm V$) 
be a basis of a free $\Z$ module $L$.
$\Atau$ is 
 a universal $\C^*$ algebra having unitary 
generators $U_i$, $1\leq i\leq n$, obeying the relations
\[
U_i U_j=e^{2\pi i \tau(\beta_i,\beta_j)} U_j U_i.
\]
It is obvious that the two definitions are equivalent.

\begin{prop}
\label{prop:taumod}
Under the conditions of Lemma \ref{lem:ration} there exists
an $N$ dimensional module $M$ over $\Atau$.
\end{prop}

{\bf Proof:} Since $N\tau$ is integral form and we have a freedom
in choosing a basis $\{\beta_i\}$ of $L$, we can choose it
so that
\[
\begin{array}{l}
N\tau(\beta_{2i-1},\beta_{2i})=q_1q_2\cdots q_i,\\[12pt] 
\tau(\beta_k,\beta_l)=0~~~{\rm unless}~~~k=2i-1~{\rm and}~l=2i~~
{\rm or}~~k=2i~{\rm and}~l=2i-1,
\end{array}
\]
where $q_1,q_2,\cdots$ are integers  (see \cite{Igusa})
and  moreover the basis $\{\beta_i\}$ respects the decomposition of $L$ into 
$L_\mu\oplus \wt{L_\mu}$.
In this basis the algebra $\Atau$ is generated by unitary
generators $U_i$ obeying the relations
\begin{equation}
\label{eq4}
U_{2i-1}U_{2i}=e^{2\pi i \frac{q_1\cdots q_i}{N}}U_{2i}U_{2i-1}
\end{equation}
( all other generators commute). Note that it may happen
that there exists an integer $m$ such that if $i>m$ then all $q_i$ are
zero. 

So we see that our algebra $\Atau$ is a tensor product of algebras
$\Atau_i$, where $\Atau_i$ is generated by two unitary generators
$U_{2i-1}$ and $U_{2i}$ obeying  relations \ref{eq4} or
$\Atau_i$ is generated by only one unitary generator (this is the case
when $i>m$, in particular if $\beta_i\in \wt{L_\mu}$).
Thus it is enough to show that we can construct finite dimensional
modules $M_i$ over $\Atau_i$ such that $(\dimm M_1)(\dimm M_2)\cdots$
divides $N$. Really, in such case we can take $M$ to be the direct sum
of $M_1\otimes M_2\otimes\cdots$ taken $\frac{N}{(\dim M_1)(\dim M_2)\cdots}$
times.

If $\Atau_i$ is generated by one unitary generator then it has a 1 dimensional
module over it, $U_i$ acts by 1. We choose $M_i$ to be this module
in this case.

 If $\Atau_i$ is generated by two unitary generators
$U_{2i-1}$ and $U_{2i}$ obeying relations \ref{eq4} then it has
a module of dimension $\frac{N}{{\rm GCD}(N,q_1\cdots q_i)}$, where
GCD stands for greatest common divisor. We choose $M_i$ to be a module of
the dimension $\frac{N}{{\rm GCD}(N,q_1\cdots q_i)}$.

Thus, it is enough to show that 
$\frac{N^m}{{\rm GCD}(N,q_1)\cdots{\rm GCD}(N,q_1\cdots q_m)}$ divides $N$,
 where $m$-the number of tori
$\Atau_i$ generated by two generators. Therefore it is enough to prove that
$\frac{{\rm GCD}(N,q_1)\cdots{\rm GCD}(N,q_1\cdots q_m)}{N^{m-1}}$ is 
an integer.

\begin{lem}
\label{lem:division}
\[
\frac{{\rm GCD}(N,q_1)\cdots{\rm GCD}(N,q_1\cdots q_m)}{N^{m-1}} 
\]
 is an integer if 
\[
\frac{q_1^2q_2}{N},\frac{q_1^3q_2^2q_3}{N^2},\cdots,
\frac{q_1^{m}q_2^{m-1}\cdots q_m}{N^{m-1}} 
\]
are integers.
\end{lem}

{\bf Proof:} Denote by $a_i$ the 
$\frac{{\rm GCD}(N,q_1\cdots q_i)}{{\rm GCD}(N,q_1\cdots q_{i-1})}$ and 
$b_i=q_i/a_i$. Then we can write
$q_1^{m}q_2^{m-1}\cdots q_m=(a_1^{m}a_2^{m-1}\cdots a_m)
(b_1^{m}b_2^{m-1}\cdots b_m)$. We will prove by induction  that
$ \frac{a_1^{k}a_2^{k-1}\cdots a_k}{N^{k-1}}$ are integers.
The initial case $k=1$ is obvious.
For $k>1$ we have
\[
\frac{a_1^{k}a_2^{k-1}\cdots a_k}{N^{k-1}}=
\frac{a_1^{k-1}a_2^{k-2}\cdots a_{k-1}}{N^{k-2}} \left(\frac{a_1\cdots a_k}{N}
\right).
\]
$\frac{a_1^{k-1}a_2^{k-2}\cdots a_{k-1}}{N^{k-2}}$ is an integer by induction
hypothesis and we also know that $\frac{N}{a_1\cdots a_k}$ are relatively prime
with $b_1,\cdots b_k$. But on the other hand
\[
\frac{q_1^{k}q_2^{k-1}\cdots q_k}{N^{k-1}}=
\frac{a_1^{k-1}a_2^{k-2}\cdots a_{k-1}}{N^{k-2}} 
\left(\frac{a_1\cdots a_k}{N}\right)
\left(b_1^{k}b_2^{k-1}\cdots b_k\right).
\]
Thus, $\frac{a_1^{k-1}a_2^{k-2}\cdots a_{k-1}}{N^{k-2}} 
\left(\frac{a_1\cdots a_k}{N}\right)$ is an integer.\\
\QED

To prove the proposition it is enough to show that
\[
\frac{q_1^2q_2}{N},\frac{q_1^3q_2^2q_3}{N^2},\cdots,
\frac{q_1^{m}q_2^{m-1}\cdots q_m}{N^{m-1}} 
\]
are all integers. Let $\{\gamma_i\}$ be a basis of $V^*$ dual 
to $\{\beta_i\}$. We know that $\mu=N e^{\iota(\tau)}\alpha$
is an integral element of $\Vev$. Denote by $k$ the dimension
of $W_\mu$ as before. 
Then, $\alpha=\pm \beta_1\wedge\beta_2\wedge\cdots\wedge\beta_k$
and the numbers
\[
\begin{array}{l}
\langle \gamma_3\wedge\cdots\wedge\gamma_k,N\iota(\tau)\alpha\rangle,
\langle \gamma_5\wedge\cdots\wedge\gamma_k,
N\frac{(\iota(\tau))^2}{2}\alpha\rangle,\\[12pt]
\cdots,
\langle \gamma_{2m+1}\wedge\cdots\wedge\gamma_k,
N\frac{(\iota(\tau))^m}{m!}\alpha\rangle
\end{array}
\]
are integers. A straightforward calculation shows that
\[
\langle \gamma_{2j+1}\wedge\cdots\wedge\gamma_k,
yN\frac{(\iota(\tau))^j}{j!}\alpha\rangle=
\pm 
\frac{q_1^{j}q_2^{j-1}\cdots q_j}{N^{j-1}}.
\]
Thus all numbers $\frac{q_1^{j}q_2^{j-1}\cdots q_j}{N^{j-1}}$ are
integers.\\
\QED

\section{Proof of Theorem \ref{thm:M1}.}
\label{sec:proof}

First, let us show that if we have a  projective
 $\At$-module with  constant curvature connection 
then the corresponding class $\mu=[E]\in K_0(\At)$ is a positive
generalized quadratic exponent.

Really, we know from formula \ref{eqrelchmu} that
\[
\mu=[E]=e^{-\iota(\theta)}ch(E)=e^{\iota(-\theta)}ch(E).
\]
Also, from formula \ref{eq:ch} we see that
\[
ch(E)=d_E e^{\frac{F}{2\pi i}}.
\]
Thus $ch(E)$ is a generalized quadratic exponent since $F$ is an element
of $\Lambda^2 V$ (recall that $V=\Lt^*$).
>From Proposition \ref{prop:closure} it immediately follows that
$\mu$ is a generalized quadratic exponent since $-\theta\in \Lambda^2 (V^*)$.
Therefore we proved that $\mu=[E]$ is a generalized quadratic exponent.
It is a positive element of $K_0(\At)$ because it represents a genuine
$\At$ module. Thus we proved the statement of Theorem \ref{thm:M1} in one
direction.

This was the easy part. The hard part is to prove the second half, that is
to show that if $\mu\in K_0(\At)$ is a positive generalized quadratic exponent
then there exists a  projective module $E$ with a 
constant curvature connection which represents the class $\mu$, {\it i.e.},
$\mu=[E]$. In the next subsection we present an explicit construction 
of such a module.

\subsection{Construction of  $\At$-module $E$ with  constant curvature
connection.}
\label{sec:constr}

In this section we will construct explicitly a $\At$-module $E$ with a
constant
curvature connection over  representing $\mu\in K_0(\At)$.
Assuming that such a module exists we see that  
$ch(E)=e^{\iota(\theta)}\mu$ is a generalized quadratic exponent
(follows from Proposition \ref{prop:closure} and the fact that
$\mu$ is a generalized quadratic exponent). Moreover,
$ch(E)_{(0)}>0$ therefore from Corollary \ref{cor:const} follows that
\[
ch(E)=e^{\iota(\theta)}\mu=d_E e^{\frac{F}{2\pi i}},
\]
where $F$ is the curvature form, and $d_E$ is the dimension of the module
$E$. Thus, reversing the previous arguments it is obvious that it is enough
to construct a  projective  $\At$- module
with  constant curvature connection  satisfying the following
properties\\
a) the curvature form is $F$;\\
b) the dimension of the module is $d_E$. \\

In Section \ref{sec:int} we defined a subspace $W_\mu\subseteq V=\Lt^*$ 
associated
with the generalized quadratic exponent $\mu$. Since $\mu\in K_0(\At)$ we
see that $\mu$ is integral. Thus, $L_\mu=L\cap W_\mu$ is the integral lattice
in $W_\mu$ by Proposition \ref{prop:intgqe}. As in section \ref{sec:int}
we denote by $k=\dimm W_\mu$ and we choose a complement $\wt{L_\mu}$ to
$L_\mu$ in $L$. Denote by $Y_\mu$ the span of $\wt{L_\mu}$ in $V$.
It is obvious that ${\Lt}^*=V=W_\mu\oplus Y_\mu$. Thus, we have a 
natural decomposition $\Lt=V^*=W^*_\mu\oplus Y^*_\mu$. Note
that the space $Y^*_\mu=U_\mu$  was defined in section \ref{sec:int}.
Since $\mu_{(k)}=d_E (e^{\frac{F}{2\pi i}})_{(k)}$
$k$ is even integer, that is $k=2p$, $p\in \Z$. Denote by 
$q$ the rank of the free abelian group $\wt{L_\mu}$. 
We have $q=n-2p$, where $n=\dimm \Lt$
the dimension of $\At$. 
Since $F^p$ is nonzero it follows that $F|_{W_\mu^*}$ 
($F$ restricted to $W_\mu^*$) is a non-degenerate 2-form.

\subsubsection{Construction of operators $\n_x$ for $x\in\Lt$.}

Let $\He$ be a Heisenberg algebra generated
by the operators $\n_x$, for $x\in W_\mu^*$ which satisfy the relation
\[
[\n_x,\n_y]=F(x,y),
\]
where $x,y\in  W_\mu^*$. The algebra  $\He$ has a unique irreducible
representation
which can be realized in the space of square integrable functions on
$\R^p$. Moreover, the action of $\n_x$ is given by an operator
\[
(\n_x(f))(z)=2\pi i \langle \phi(x), z\rangle f(z) + 
\sum_i \psi_i(x)\frac{\partial f(z)}{\partial z_i},
\]
where $\phi: W_\mu\rightarrow (\R^p)^* $ is some linear map, and
$\psi_i: W_\mu\rightarrow \R$ are some linear functions on $W_\mu^*$.
In particular, we see that these operators preserve the space
of Schwartz functions on $\R^p$. 

The above construction provides us with the action of the operators
$\n_x$ for $x\in W_\mu^*$ only. First we will extend the above construction 
to obtain an action
of all operators $\n_x$, $x\in \Lt=W_\mu^*\oplus Y_\mu^*$. Then, we 
will obtain an action of some non-commutative torus $\Att$ so that
$\n$ becomes an $\Att$ connection.

We extend the space from the space of Schwartz functions on $\R^p$
to the space of Schwartz functions on 
$\R^p\times \wt{L_\mu}=\R^p\times \Z^q$. Denote it by $H$.
If $x\in \Lt=W_\mu^*\oplus Y_\mu^*$ we denote by $x_W$ the projection
of $x$ on $W^*_\mu$ and by $x_Y$ the projection of $x$ on $Y^*_\mu$
(obviously $x=x_W+x_Y$).
We define the action of $\n_x$ on an element $f(z,a)\in H$, where $z\in\R^p$ 
and $a\in \wt{L_\mu}$, as follows
\begin{equation}
\label{eq:actn}
(\n_x(f))(z,a)=(\n_{x_W}(f))(z,a)+2\pi i \langle x_Y,a\rangle f(z,a),
\end{equation}
where the action of $\n_{x_W}$ is the same as above (only along $z$'s).
Notice that the operators $\n_x$, $x\in\Lt$ satisfy the commutation relations
\[
[\n_x,\n_y]=[\n_{x_W},\n_{y_W}]=F(x_W,y_W)=F(x,y),~~~x,y\in\Lt,
\]
since $\n_{x_Y}$ (recall 
that $(\n_{x_Y}(f))(z,a)=2\pi i \langle x_Y,a\rangle f(z,a)$) 
commutes with $\n_{y}$ for any $y\in\Lt$. 

Thus, we constructed operators
$\n_x$, $x\in\Lt$ which satisfy the desired commutation relations.
Next, we will construct operators acting on the space $H$ which 
generate  a non-commutative torus $\Att$ such that
\begin{equation}
\label{cond}
\begin{array}{l}
{\rm \bf first}~~~\n~~{\rm is~an}~~\Att~~{\rm connection}\\ 
{\rm \bf second}~~~\theta-\wt{\theta}~~{\rm is~an~element~of}~~
\Lambda^2 W^*_\mu.
\end{array}
\end{equation}

\subsubsection{Construction of operators satisfying  conditions
\ref{cond}.}

Let us  choose a basis $\beta_1,\cdots,\beta_{2p}$ of $L_\mu$ and
a basis
$\beta_{2p+1},\cdots,\beta_{2p+q}$ of $\wt{L_\mu}$.
We will construct operators $V_i$, $1\leq i\leq 2p+q$
acting on $H$ which generate a non-commutative torus $\Att$ which 
satisfies conditions \ref{cond}.

\begin{lem}
\label{lem:oper1}
For $1\leq i\leq 2p$ there exists an operator ${\wt{V_i}}$ acting on $H$
such that\\  
$
{\rm first}~~~({\wt{V_i}}(f))(z,a)=
e^{2\pi i \chi_i(z)}f(z+y_i,a),~{\it where}~z\in \R^p,~~
a\in \wt{L_\mu},
$\\   
\hspace*{0.5in} 
for some $y_i\in \R^p$ and some linear function $\chi_i\in(\R^p)^*$;\\
{\rm second} 
$[\n_x,{\wt{V_i}}]=2\pi i \langle x,\beta_i \rangle {\wt{V_i}}$, for
$x\in \Lt$.
\end{lem}

{\bf Proof:} 
Let us introduce an operator $W(y,\chi)$, where 
$y\in \R^p$ and $\chi\in (\R^p)^*$
\[
(W(y,\chi)f)(z,a)=e^{2\pi i \chi(z)}f(z+y,a).
\]
A straightforward calculation shows  that
$[\n_x,W(y,\chi)]W(y,\chi)^{-1}$ is an operator of multiplication by a 
real number and moreover we obtain a non-degenerate pairing
between the spaces $W_\mu^*$ and $\R^p\oplus(\R^p)^*$. Thus an choosing 
appropriate element $y_i\in \R^p$ and $\chi_i\in (\R^p)^*$
we can put ${\wt{V_i}}=W(y_i,\chi_i)$.\\
\QED

We define the operators $\wt{V_i}$ for $1\leq i\leq 2p$ as in the above lemma.
We  define the operators $\wt{V_i}$ for $2p+1\leq i\leq 2p+q$ acting on $H$
 by the formula
\[
(\wt{V_i}f)(z,a)=f(z,a-\beta_i).
\]

\begin{lem}
\label{oper2}
For any $1\leq i \leq n=2p+q$, and any $x\in\Lt$ we have
\begin{equation}
\label{eq:comrel}
[\n_x,\wt{V_i}]=2\pi i \langle x,\beta_i\rangle \wt{V_i}.
\end{equation}
\end{lem}

{\bf Proof:} For $i\leq 2p$ formula \ref{eq:comrel} follows from
Lemma \ref{lem:oper1}. For $i>2p$ formula \ref{eq:comrel} follows
from an easy straightforward calculation.\\
\QED

It is easy to check that the operators $\wt{V_i}$ are generators
of some non-commutative torus. Moreover, these operators satisfy 
the first condition in (\ref{cond}) but they do not satisfy the
second condition. To remedy this we will modify the operators $\wt{V_i}$
replacing them with operators
$V_i=e^{2\pi i l_i(\cdot)}\wt{V_i}$.

If $l\in Y_\mu^*$ then the operator $e^{2\pi i l(\cdot)}\wt{V_i}$ acts
on $H$ by the formula
\[
(e^{2\pi i l(\cdot)}\wt{V_i}(f))(z,a)=e^{2\pi i l(a)}(\wt{V_i}(f))(z,a).
\]
Moreover, we have
\[
[\n_x,e^{2\pi i l(\cdot)}\wt{V_i}]= 2\pi i \langle x,\beta_i\rangle
e^{2\pi i l(\cdot)}\wt{V_i}
\]
which follows from an easy straightforward calculation (since the operator
$e^{2\pi i l(\cdot)}$ commutes with the operators $\n_x$, $x\in\Lt$).

\begin{prop}
\label{mainconst}
For $1\leq i\leq 2p+q$ there exists a linear function 
$l_i\in Y_\mu^*$ on $Y_\mu$ such that\\
if we define 
$V_i=e^{2\pi i l_i(\cdot)}\wt{V_i}$
 then
\begin{equation}
\label{eq:tettild}
V_i V_j=e^{2\pi i \wt{\theta}_{ij}}V_j V_i,
\end{equation}
and $\theta-\wt{\theta}$ is an element of $\Lambda^2 W^*_\mu$
\end{prop}

{\bf Proof:} First, it is easy to see that there exists a 2-form
$\sigma\in \Lambda^2 \Lt$ such that 
\[
\wt{V_i} \wt{V_j}=e^{2\pi i \sigma_{ij}}\wt{V_j} \wt{V_i}.
\]
An easy calculation shows that the operator $e^{2\pi i l(\cdot)}$
commutes with the operators $\wt{V_i}$ for $i\leq 2p$. If $i>2p$ then
we have
\[
(\wt{V_i}\circ e^{2\pi i l(\cdot)})= e^{-2\pi i l(\beta_i)}
(e^{2\pi i l(\cdot)}\circ\wt{V_i}).
\]
This gives us that\\
 if $i,j\leq 2p$ then
\begin{equation}
\label{eq:ijleqtp}
V_i V_j=e^{2\pi i \sigma_{ij}}V_j V_i;
\end{equation}
if $i\leq 2p$ and $j>2p$ then
\begin{equation}
\label{eq:ileqtpjgtp}
V_i V_j=e^{2\pi i (\sigma_{ij}+l_i(\beta_j))}V_j V_i;
\end{equation}
and if $i,j>2p$ then
\begin{equation}
\label{eq:ijgtp}
V_i V_j=e^{2\pi i (\sigma_{ij}+l_i(\beta_j)-l_j(\beta_i))}V_j V_i.
\end{equation}
For $1\leq i\leq 2p$ we define $l_i\in Y_\mu^*$ by the formula
\[
l_i(\beta_j)=\theta(\beta_i,\beta_j)-\sigma_{ij}
\]
on the basis $\{\beta_j\}$, $j>2p$ of $Y_\mu$.
For $2p<i\leq2p+q$ we define $l_i\in Y_\mu^*$ by the formula
\[
l_i(\beta_j)=\frac{1}{2}(\theta(\beta_i,\beta_j)-\sigma_{ij})
\]
on the basis $\{\beta_j\}$, $j>2p$ of $Y_\mu$.

Equations \ref{eq:ijleqtp}, \ref{eq:ileqtpjgtp}, and \ref{eq:ijgtp}
show that
\[
V_i V_j=e^{2\pi i \theta(\beta_i,\beta_j)}V_j V_i
\]
if either $i$ or $j$ greater then $2p$ and
\[
V_i V_j=e^{2\pi i \sigma_{ij}}V_j V_i
\]
if $i,j\leq 2p$. 
Thus we constructed the linear functions $l_i\in Y_\mu^*$ such that
the conditions \ref{eq:tettild} are satisfied.\\
\QED

We define the operators $V_i$ as in the above Lemma. We easily
see that the operators $V_i$ generate a non-commutative torus $\Att$,
where $\wt{\theta}(\beta_i,\beta_j)=\sigma_{ij}$ if both $i,j\leq 2p$
and $\wt{\theta}(\beta_i,\beta_j)={\theta}(\beta_i,\beta_j)$ otherwise.

Thus, the operators $V_i$ satisfy the condition \ref{cond}.

\subsubsection{Construction of a  projective $\Att$
module.}

Now we will identify our construction with the construction given in
\cite{Rf1}. 

Let $G$ be a central extension of the abelian group
$\R^p\times \wt{L_\mu}\times (\R^p)^*\times (\wt{Y_\mu}/L_\mu^*)$
given by the natural pairing between $\R^p\times \wt{L_\mu}$ and
$(\R^p)^*\times (\wt{Y_\mu}/L_\mu^*)$. We see that $G$ is a
Heisenberg group and it acts naturally on $H$. 
We denote this representation by
$\rho$. 
Moreover, for each $V_i$ there 
exists a unique element $g_i\in G$ such that $\rho(g_i)=V_i$. 
One can easily recognize  the construction of elementary modules
over non-commutative tori in M. Rieffel's paper \cite{Rf1}. 

Thus, 
choosing an appropriate space of functions on $\R^p\times \wt{L_\mu}$
we get a  projective  $\Att$module $\wt{E}$ with  constant
curvature connection $\n$ 
 such that\\
{\it first}, the curvature of $\n$ is $F$;\\
{\it second}, $\theta-\wt{\theta}$ is an element of $\Lambda^2 W_\mu^*$.

Next, we would like to find explicitly the class $[\wt{E}]$ in $K_0(\Att)$.
Note, that in our construction of module $\wt{E}$ we canonically 
identified the space $\Lt$ with the space $L_{\wt{\theta}}$. Thus, we can
think about $[\wt{E}]$ as an integral element of 
$\Lambda^{even}\Lt^*=\Lambda^{even}V$. To find the
class $[\wt{E}]$ we would have to do some calculations. Fortunately,
they were already done by M. Rieffel in \cite{Rf1}. So, we will apply his
results to our case. 

Let us remind that in paper \cite{Rf1} M. Rieffel introduced a linear
map $\wt{T}:L_{\wt{\theta}}\rightarrow \R^p\times \R^q\times(\R^p)^*$.
In our notation $L_{\wt{\theta}}$ is canonically identified with
$\Lt=V=W_\mu\oplus Y_\mu$ and $\R^q$ with $Y_\mu$. Thus, in our terms
we have a linear map 
$\wt{T}:W_\mu\oplus Y_\mu\rightarrow \R^p\times Y_\mu\times(\R^p)^*
=\R^p\times(\R^p)^*\times Y_\mu$. It is easy to see from the 
explicit construction of the operators $V_i$
that $\wt{T}$ maps $W_\mu$ to $\R^p\times(\R^p)^*$,
and $Y_\mu$ to $Y_\mu$. Moreover  the restriction
of $\wt{T}$ on $Y_\mu$ is the identity map. 

M. Rieffel found in \cite{Rf1} that 
\begin{equation}
\label{eq:class}
[\wt{E}]=d\prod_{j=1}^p\bar{Y}_j\wedge\bar{Y}_{j+p},
\end{equation}
where $d={\rm det}(\wt{T})$ and
\begin{equation}
\label{eq:Y}
\bar{Y}_j=
\begin{array}{l}
\wt{T}^{-1}(\bar{e}_j)~~~{\rm for}~~~1\leq j\leq p\\[10pt]
\wt{T}^{-1}(e_{j-p})~~~{\rm for}~~~p+1\leq j\leq 2p,
\end{array}
\end{equation}
where $\{e_j\}$ is a basis of $\R^p$ and $\{\bar{e}_j\}$ is the dual basis
of $(\R^p)^*$.

Since $\wt{T}$ is the identity map on $Y_\mu$ we see that
\[
{\rm det}(\wt{T}^{-1})=\pm\frac{\prod_{j=1}^p\bar{Y}_j\wedge\bar{Y}_{j+p}}
{\alpha},
\]
where $\alpha$ is the volume form on $W_\mu$ (see Proposition \ref{prop:gqe}
for the definition of $\alpha$). Note, we put a $\pm$ sign because we
do not want to specify precisely how to pick a volume form. The lattice
$L_\mu$ specifies the volume form upto a sign. Later it will be easy to 
make the right choice of the sign so that everything would agree with
M. Rieffel's paper \cite{Rf1}.
We get
\[
d={\rm det}(\wt{T})=1/{\rm det}(\wt{T}^{-1})= \pm\frac{\alpha}
{\prod_{j=1}^p\bar{Y}_j\wedge\bar{Y}_{j+p}}.
\]
Thus we obtain that
\begin{equation}
\label{eq:classEt}
[\wt{E}]=\pm\alpha.
\end{equation}

\subsubsection{Construction of a  projective $\At$
module $E$.}

{}From the above results and Proposition \ref{prop:intgqe} we see that
if we make the right choice of the sign (so that $[\wt{E}]=\alpha$) then 
\begin{equation}
\label{defN}
N=\frac{\mu_{(2p)}}{\alpha}
\end{equation}
is a positive integer. Moreover, we have
\[
d_{\wt {E}} e^{\frac{F}{2\pi i}}=e^{\iota(\theta)}\mu
\]
and
\[
d e^{\frac{F}{2\pi i}}=e^{\iota(\wt{\theta})}\alpha.
\]
Therefore, 
\begin{equation}
\label{equal}
e^{\iota(\theta)}\mu=\left(\frac{d_{\wt
E}}{d}\right)e^{\iota(\wt{\theta})}\alpha.
\end{equation}
{}From equations \ref{equal} and \ref{defN} we see that
$\frac{d_{\wt E}}{d}=N$.

One can easily check that the conditions of Lemma \ref{lem:ration}
are satisfied. Therefore  
$N(\theta-\wt{\theta})$ is an integral element of $(W_\mu)^*=X_\mu$.
{}From Proposition \ref{prop:taumod} follows that there exists
$N$-dimensional
module $M$ over ${\cal A}_{\theta-\wt{\theta}}$.
Denote by
\begin{equation}
\label{E}
E=\wt{E}\otimes M.
\end{equation}
{}From Proposition 5.4 and Theorem 5.6 in M. Rieffel's paper 
\cite{Rf1} it follows
that $E$ is a  projective module over $\At$ (since
$\theta=\wt{\theta}+(\theta-\wt{\theta})$) with constant curvature 
connection with the curvature
 given by  formula 
\[
\Omega=F\otimes {\rm Id}_M=F
\]
 and
\[
{\rm ch}(E)=\dim(M){\rm ch}(\wt{E}).
\]
Therefore we see that
\[
{\rm ch}(E)=N{\rm ch}(\wt{E})=
N e^{\iota(\wt{\theta})}\alpha=
N \left(\frac{d}{d_E}\right)e^{\iota(\theta)}\mu=
e^{\iota(\theta)}\mu.
\]
Thus we constructed a  projective  $\At$-module $E$
with  constant curvature connection 
 such that
$[E]=\mu$. This finishes the proof of Theorem \ref{thm:M1}\\
\QED

{\bf Acknowledgements.} We would like to express our deep gratitude to 
Dmitry Fuchs, Anatoly Konechny and Marc Rieffel for
numerous fruitfull discussions.

\end{document}